\documentclass{article}

\usepackage{authblk}
\usepackage{amsmath,amssymb,amsthm}
\usepackage[usenames,dvipsnames]{color}
\usepackage{times}
\usepackage{amsxtra}
\usepackage{enumerate}
\usepackage{pdfsync}
\usepackage{enumitem}

\usepackage{graphicx}
\usepackage{wrapfig}
\usepackage{caption,subcaption}

\usepackage{tabularx}
\usepackage{caption}
\usepackage{float}
\usepackage[T1]{fontenc}

\newtheorem{asm}{Assumption}
\newtheorem{rmk}{Remark}
\newtheorem{thm}{Theorem}
\newtheorem{lemma}{Lemma}

\newtheorem{propos}{Proposition}
\newtheorem{clr}{Corollary}

\numberwithin{equation}{section}
\numberwithin{thm}{section}
\numberwithin{rmk}{section}
\numberwithin{lemma}{section}
\numberwithin{clr}{section}
\numberwithin{Def}{section}
\newcommand{\PP}{\mathbb{P}}
\newcommand{\EE}{\mathbb{E}}

\newcommand{\RR}{\mathbb{R}}
\newcommand{\eps}{\varepsilon}
\newcommand{\td}{\tau_d}
\newcommand{\la}{\lambda}
\usepackage[natbibapa,nodoi]{apacite}
\begin{document}
	
	\title{On the optimally controlled stochastic shallow lake}


    %
	
	\author[a]{Angeliki Koutsimpela}
	\author[b]{Michail Loulakis}
		\affil[a,b]{School of Applied Mathematical and Physical Sciences, National Technical University of Athens, 15780 Athens, Greece}
		\affil[b]{Institute of Applied and Computational Mathematics, Foundation of Research and Technology Hellas,
			 70013 Heraklion Crete, Greece}	
	
	\maketitle
	
	\begin{abstract}
		We consider the stochastic control problem of the {\em shallow lake} and continue the work of \cite{KLS} in three directions. First, we generalise the characterisation of the value function as the viscosity solution of a well-posed problem to include more general recycling rates. Then, we prove approximate optimality under bounded controls and we establish quantitative estimates. Finally, we implement a convergent and stable numerical scheme for the computation of the value function to investigate properties of the optimally controlled stochastic shallow lake. This approach permits to derive tails asymptotics for the invariant distribution and to extend results of \cite{GKW} beyond the small noise limit.\\
  
       
	\end{abstract}
 


	\section{Introduction}
	
The shallow lake problem is a well-known problem of the environmental  economy with a great mathematical interest. Pollution of shallow lakes is caused by human activity, e.g. the use of fertilisers and the increased inflow of waste water from industries and human settlements, and is usually quantified by the concentration of phosphorus.   The amount of phosphorus in algae, $x(t)$, is usually modelled by the non-linear stochastic differential equation:
\begin{equation} \label{sldyn}
	\begin{cases} dx(t)= \left( u(t)-b x(t)+r\big(x(t)\big)\right) dt+\sigma x(t)dW_t, & \\
		x(0) = x\ge 0. & \end{cases}
\end{equation}
The first term, $u:\; [0,\infty)\rightarrow (0,\infty)$, in the drift part of the dynamics, represents the exterior load of phosphorus as a result of human activities. The second term is the rate of loss $b x(t)$, which is due to sedimentation, outflow and sequestration in other biomass. The third term, $r\big(x(t)\big)$, is the rate of recycling of phosphorus on the bed of the lake. This term is assumed to be a sigmoid function (see \cite{carp99}) and the typical choice in the literature is the function $x\mapsto \frac{x^2}{x^2+1}$. An uncertainty in the rate of loss is inserted in the model through a linear multiplicative gaussian white noise with  intensity $\sigma$.

The economics of the lake arise from its conflicting services to the community. On the one hand, a clear lake is an ecological resource of recreational activities. On the other hand, the lake can serve as a sink for agricultural and industrial waste. When the users of the lake cooperate, the loading strategy,  $u$, can be used as a control to maximise the benefit from the lake. Assuming an infinite horizon, this benefit is typically chosen as
\begin{equation}\label{Jxu}
	J(x;u)=\EE_x \left[\int_0^\infty e^{-\rho t}\big(\ln u(t)-cx^2(t)\big)\, dt \right],
\end{equation}
where $\rho >0$ is the discount rate and $x(\cdot)$ is the solution to (\ref{sldyn}), for a given exterior loading (control)  $u(\cdot)$, and initial state, $x\ge 0$. The total benefit of the lake increases with the loading of phosphorus as $\ln u(t)$, but at the same time decreases with the existing amount of phosphorus as $-cx^2(t)$, due to implied deterioration of its ecological services. The positive parameter $c$ reflects the relative weight of this component.

For the optimal management of the lake when the initial state is $x$, we need to maximise the total benefit over a set of admissible controls, $\mathfrak{U}_x$. Thus, the value function of the problem is  
\begin{equation}\label{ihvf}
	V(x) = \underset{u\in \mathfrak{U}_x} {\sup}\ J(x;u).
\end{equation}
\vskip.05in

Therefore, the shallow lake problem becomes a problem of control theory or a differential game in the case where we have competitive users of the lake (\cite{carp99, Brock, Xep}).

The deterministic ($\sigma=0$) version of the problem has been extensively studied, not only as an interesting problem in environmental economics but also as a prototype example where a sigmoid term in the dynamics may lead to the existence of multiple equilibria and associated basins of attraction. Depending on the parameters, the shallow lake problem may have two different equilibria and a Skiba point, i.e., an initial state for which distinct optimal solutions to \eqref{ihvf} exist. The leftmost (oligotrophic) equilibrium point of the system of the lake corresponds to a lake with low concentration of phosphorus, while the rightmost one (eutrophic) corresponds to a lake with high concentration of phosphorus. At the Skiba point, there are two different optimal strategies, each one driving the system to a different equilibrium and the value function is not differentiable thereat. 

The range of parameters for which Skiba points appear has been explored in  \cite{W03} and \cite{WBif}. Properties of the value function of the deterministic shallow lake problem have been proved in  \cite{KosZoh}. The existence of optimal control is usually taken as a hypothesis in the literature and the optimal dynamics of the lake is studied mostly through the necessary conditions determined by the Pontryagin maximum principle, and the equilibrium points of the corresponding dynamical system  (see \cite{W03, Xep}).  A rigorous answer to this question was given in \cite{Bar20} and \cite{Bar21}, albeit under restrictions that do not fully cover the range of  the parameters for which Skiba points are present. 
 
Over the last decade, there has been increasing interest in the stochastic ($\sigma\neq 0$) version of the problem.  Deterministic systems with two equilibrium points and one Skiba point have a fundamentally different behaviour  from their stochastic counterparts. Specifically, random fluctuations drive the stochastic system from one equilibrium point to the other (metastability).  In the context of the shallow lake, this phenomenon is studied numerically in \cite{GKW}, where the value function in \eqref{ihvf} is approximated for small $\sigma$ based on heuristic methods of perturbation analysis. In \cite{KLS} the authors characterise the value function of the stochastic shallow lake problem as the unique (in a suitable class) state-constraint viscosity solution of the Hamilton-Jacobi-Bellman (HJB) equation
\begin{equation}\label{OHJB}
	\rho V-\left(r(x)-bx\right)V_x+\left( \ln(-V_x)+cx^{2}+1\right)-\dfrac{1}{2}\sigma^{2}x^{2}V_{xx}=0,
\end{equation}
and analytically derive properties of the value function in the case $r(x)=x^2/(1+x^2)$.

The present work continues the work of  \cite{KLS}. First, it extends the results therein to include much more general recycling rates as made precise in Assumption \ref{r}, as well as the penalty parameter $c$ in \eqref{Jxu}, which cannot be scaled away with a suitable change of variables.  Then, it adds some new analytical results regarding the approximate optimality of bounded controls (Lemma \ref{Napprox}) and the tail behaviour of the invariant distribution of the optimally controlled stochastic lake (Proposition \ref{asympt}). Finally, we use a convergent, monotone, Barles-Souganidis scheme to compute the value function of \eqref{ihvf} as the relevant viscosity solution of \eqref{OHJB}. This approach relies heavily on the aforementioned rigorous results and because it does not invoke perturbation expansions, it permits to numerically investigate properties of the optimally controlled stochastic shallow lake beyond the small noise regime. As an example, we present the effect of noise intensity, $\sigma$, of the penalty parameter $c$, and of the discount rate, $\rho$, on the number and location of modes of the invariant distribution. We also present typical paths of the optimally controlled stochastic shallow lake showing the transitions between oligo/eu-trophic states, as well as statistics of the transition times. Analytical results are collected in the following section and numerical results are presented in Section \ref{numerics} of this article.
		
	\section{Analytical Results}
In this section, we introduce the details of the model, some necessary notation, and we present the rigorous part of this work. The section contains two kinds of results: generalisations of those in \cite{KLS} and new ones. Regarding the proofs of the former, naturally, some arguments do not depend on the precise form of the recycling rate and go through verbatim simply by substituting $x^2/(1+x^2)$ by $r(x)$, while others need to be modified to a lesser or greater extent. For the sake of completeness, we transcribe here most of the results in \cite{KLS} in their generalised form but, to avoid cumbersome repetitions, we only present those proofs that differ from the original ones. To keep the flow of results uninterrupted, statements are collected in  subsection \ref{statement}, while all proofs provided are collected in subsection \ref{proofs}.
	
We assume that there exists a filtered probability space $( \Omega, \mathcal{F},\{ \mathcal{F}_t\}_{t\ge 0},\PP )$
satisfying the usual conditions, and a Brownian motion $\{W_t: t\ge 0\}$ defined on that space. 
An admissible control   $u(\cdot)\in\mathfrak{U}_x$ is an $\mathcal{F}_t$-adapted, $\PP$-a.s. locally integrable process with values in $U=(0,\infty)$, satisfying 
\begin{equation}\label{ac_constraint}
	\mathbb{E} \left[ \int_{0}^{\infty}e^{-\rho t}\ln u(t)dt \right] < \infty,
\end{equation}
such that  the problem (\ref{sldyn}) has a unique strong solution $x(\cdot)$. We will also find useful the positive processes 
\begin{equation}\label{ZM}
		Z_t=e^{\sigma W_t-\left(b+\sigma^2/2\right)t} \quad \text{ and } \quad M_t(u)=\int \limits_{0}^{t}\frac{Z_t}{Z_s}u(s)ds.
\end{equation}
\vskip.05in
Regarding the recycling rate, $r$, the following assumptions are made throughout the paper. Note that they are satisfied by most sigmoid functions that are typically used in applications.
\begin{asm}\label{r}
	The rate of recycling $r(x)$ satisfies the following:
	\begin{enumerate}
		\item $r$ is locally Lipschitz and nondecreasing 
		\item $r(0)=0$ and $r(x)<(b+\rho)x$ close to $0$
		\item $a:=\lim \limits_{x\rightarrow \infty}r(x)<\infty$
		\item The limit $\lim \limits_{x\rightarrow \infty}(a-r(x))x=:C$ exists and is a finite, necessarily nonnegative, real number.
	\end{enumerate}
\end{asm}
\subsection{Statements}\label{statement}
Proposition \ref{propdyn} collects some useful facts about the set of admissible controls and path properties of the controlled system \eqref{sldyn}. Precisely, the set of admissible controls does not depend on the initial state, while paths of \eqref{sldyn} remain nonnegative and satisfy a comparison principle.
\noindent
\begin{propos}\label{propdyn}
	\begin{enumerate}[before=\leavevmode,label=\upshape(\roman*),ref=\thepropos (\roman*)]
		\item \label{positive} If $x\geq 0$, $u\in \mathfrak{U}_x,$ and $x(\cdot)$ is the solution to  (\ref{sldyn}), then
		$
		\PP\big[x(t)\ge 0,\ \forall t\ge 0\big]=1.
		$ In particular, 	$
		\PP\big[x(t)\ge M_t(u),\ \forall t\ge 0\big]=1.
		$
		\item \label{uxisu}	For all $x,y\ge 0$,  $\mathfrak{U}_x=\mathfrak{U}_y=:\mathfrak{U}.$ 
		\item \label{monotone}
		Suppose $x(\cdot),\ y(\cdot)$ satisfy (\ref{sldyn}) with controls $u_1,\ u_2\in\mathfrak{U}$, respectively, and $x(0)=x$, $y(0)=y$.
		If $x\le y$ and $\PP\big[u_1(t)\le u_2(t),\ \forall t\ge 0\big]=1$, then
		\[
		\PP\big[y(t)-x(t)\ge (y-x)Z_t,\ \forall t\ge 0\big]=1.
		\]
	\end{enumerate}
\end{propos}

\noindent
Propositions \ref{propv}, \ref{vprop}, and \ref{st_eq} describe properties of the value function, $V$. Note that these properties are derived directly from the definition of $V$ in \eqref{ihvf} by means of stochastic analysis, so they are not a consequence of any differential equation, such as \eqref{OHJB}, that $V$ may satisfy. On the contrary, they are used as a crucial input in the characterisation of the value function as a viscosity solution to  \eqref{OHJB}, as they ensure that the associated Hamiltonian of the control problem is finite, and they outline a class of functions among which there is uniqueness of solutions to \eqref{OHJB}, thus singling out the relevant one. Notice also that Propositions \ref{propv} and \ref{vprop} do not require $\sigma>0$, so they can be used in the deterministic problem, as well. Let
\begin{equation}\label{K}
	A:=\frac{c}{\rho+2b-\sigma^2}.
\end{equation}
\begin{propos}\label{propv}
		\begin{enumerate}[before=\leavevmode,label=\upshape(\roman*),ref=\thepropos (\roman*)]
		\item \label{restr1}	$V(x)> -\infty$ if and only if $\sigma^2< \rho+2b$. 
		\item \label{Vdec}	The function $x\mapsto V(x)+Ax^2$, where $A$ is defined in \eqref{K}, is decreasing on $[0,+\infty)$.
		\item \label{V0fin}
		The value function at zero satisfies
		$
		V(0) \le\frac{1}{\rho}\ln\left(\frac{b+\rho}{\sqrt{2ec}}\right).
		$ 
		\item \label{specialdpp} 
		Fix $x_1,x_2\in[0,\infty)$ with $x_1<x_2$, and, for  $u\in\mathfrak{U}$, let $x(\cdot)$ be the solution to (\ref{sldyn}) with control $u$ and $x(0)=x_1$.  If $\tau_u$ is 
		the hitting time of $x(\cdot)$ on $[x_2,+\infty)$, that is,
		$
		\tau_u=\inf\{t\ge 0: x(t)\ge x_2\},
		$
		then
		\begin{equation}\label{babydpp}
			V(x_1)=\sup_{u\in\mathfrak{U}}\EE\left[\int_0^{\tau_u} e^{-\rho t}\big(\ln u(t)-cx^2(t)\big)\, dt+e^{-\rho\tau_u}V(x_2)\right].
		\end{equation}
	\end{enumerate}
\end{propos}

\begin{propos}\label{vprop}
	Suppose $0\leq \sigma^2<\rho+2b$.
	\begin{enumerate}[label=\upshape(\roman*),ref=\thepropos (\roman*)]
	\item \label{vprop1} There exist constants $K_1, K_2>0$, such that, for any $x\ge 0$, we have
	\begin{equation}\label{Vbounds}
		K_1\ \le\ V(x)+A\left(x+\frac{a}{b+\rho}\right)^2+\frac{1}{\rho}\ln \left(x+\frac{a}{b+\rho}\right)\ \le\ K_2.
	\end{equation}
	\item \label{vprop2} There exist a constant  $C_1>0$ and a function $c: \; [0,+\infty)\rightarrow (0,\infty)$ with $\lim \limits_{x\rightarrow 0}c(x)=e^{-(\rho V(0)+1)}$  such that, for any $x_1,x_2\in [0,+\infty)$ with $x_1<x_2$,
	\begin{equation}\label{DVbounds}
		\frac{V(x_2)-V(x_1)}{x_2-x_1}\le -c(x_2) \le -C_1.
	\end{equation}
\end{enumerate}
\end{propos}

\begin{propos}\label{st_eq}
	Suppose $0< \sigma^2<\rho+2b$.\\
	There exists an increasing function $L_{\sigma}: [0,\infty) \rightarrow \RR $ with $\lim \limits_{x\rightarrow 0} L_{\sigma}(x)=e^{-(\rho V(0)+1)}$ such that, for any $x_1,x_2\in [0,\infty)$ with $x_1<x_2,$
	\begin{equation}
		\frac{V(x_2)-V(x_1)}{x_2-x_1}\geq -L_{\sigma}(x_2)
	\end{equation}
\end{propos}
An immediate consequence of Propositions \ref{vprop}(ii) and \ref{st_eq} is the following corollary.
\begin{clr}\label{coro1}
If \ $0< \sigma^2<\rho+2b$, $V$ is differentiable at zero and 
	\begin{equation}
		\ln\big(-V'(0)\big)+\rho V(0)+1=0.
		\label{bc}
	\end{equation} 
\end{clr}
\begin{rmk}
Equation \eqref{bc} can be perceived as a nonlinear mixed boundary condition satisfied by the value function. It remains true even when $\sigma=0$, but the proof is more involved and is part of a forthcoming work on the deterministic shallow lake problem. A discretised version of it is used in the numerical scheme of Section \ref{numerics}.
\end{rmk}
\noindent
Theorem \ref{constrained_vs} states that the value function of the stochastic shallow lake problem is a constrained viscosity solution of the Hamilton-Jacobi-Bellman (HJB) equation \eqref{OHJB} on $[0,+\infty)$, i.e. a viscosity subsolution on $[0,+\infty)$ and a viscosity supersolution on $(0,+\infty)$. Theorem \ref{comparison} is a comparison principle that guarantees uniqueness of such solutions in a suitable class of functions. In view of Proposition \ref{vprop} and Corollary \ref{coro1}, Theorems \ref{constrained_vs} and \ref{comparison} characterise the value function as the unique constrained viscosity solution of \eqref{OHJB} on $[0,+\infty)$ that satisfies all the conditions of Theorem \ref{comparison}.

\begin{thm}\label{constrained_vs}
	If $0<\sigma^2 < \rho +2b$, the value function $V$  is a continuous in $[0,\infty) $ constrained viscosity solution of the  equation \eqref{OHJB} in  $[0,\infty)$. 
\end{thm}
\begin{thm}\label{comparison} 	Let $0<\sigma^2 < \rho +2b$ and assume that  $u\in C([0,\infty))$ is a bounded from above  strictly decreasing subsolution of \eqref{OHJB} in $[0, \infty)$ and $v\in C([0,\infty))$ is a bounded from above strictly decreasing supersolution of
	\eqref{OHJB} in $(0, \infty)$ such that  $v \geq -c_1(1+x^2)$ and  $Du\leq -\frac{1}{c_2}$
	in the viscosity sense,
	for $c_1$, $c_2$  positive constants. Then
	$u \leq v$ in $[0, \infty)$. 
\end{thm}
\noindent
An important ingredient for the proof of Theorem \ref{comparison} is the following lemma \ref{u_v}.
\begin{lemma} \label{u_v} Suppose $0\leq \sigma^2 <\rho+2b$ and $u$, $v$ satisfy the assumptions of Theorem \ref{comparison}.
	Then $\psi=u-v$ is a subsolution of 
	\begin{equation}\label{u-v}
		\rho \psi +bxD\psi-\left(a+c^*\right)|D\psi| -\frac{1}{2}\sigma^2x^2 D^2 \psi = 0 \,\, \mbox{ in } \,\, [0, \infty).
	\end{equation}
\end{lemma}
\begin{rmk}
Theorems \ref{constrained_vs} and \ref{comparison} remain valid even when $\sigma=0$. In fact, uniqueness can be established in a slightly larger class. However, the proofs in the deterministic case are quite different and, like Corollary \ref{coro1}, they are part of a forthcoming work.
\end{rmk}
\noindent
The stability property of viscosity solutions yields the following corollary. 

\begin{clr}\label{asymptotic_beh}  As  $\sigma \rightarrow 0$, the value function $V=V_{\sigma}$ defined by (\ref{ihvf}) converges locally uniformly to the constrained viscosity solution $V_0$   of the deterministic  shallow lake equation in $[0, \infty)$,
	$$\rho V_0=\left( r(x)-bx\right)V'_0-\Big( \ln(-V'_0)+cx^{2}+1\Big).$$
\end{clr}

\noindent
Finally, the next theorem describes the exact asymptotic behaviour of $V$ at $+\infty$. 			In Section 4, we present and implement a monotone numerical scheme approximating (\ref{ihvf}). 
Relation (\ref{asympto_behav}) is crucial  for the accurate computation of $V$ in this setting, because it suggests the boundary condition at the right end of the computational domain.

\begin{thm}\label{asymptotic}
	As $x\rightarrow \infty$,
	\begin{equation}\label{asympto_behav}
		V(x)= -A\left(x+\frac{a}{b+\rho}\right)^2-\frac{1}{\rho}\ln \left[2A(x+\frac{a}{b+\rho})\right]+K+o(1). 
	\end{equation}
where
\begin{equation}\label{K1}
	K=\frac{1}{\rho}\left(\frac{2b+\sigma^2}{2\rho}-\frac{Aa^2(\rho+2b)}{(b+\rho)^2} -1+2AC\right) 
\end{equation}
\end{thm}

\noindent
We now proceed to the statement of two new results. 
Lemma \ref{Napprox} states that we may essentially achieve the optimal total benefit using a bounded control $u$. For $N>0$, define the set of $N-$bounded admissible controls
$$
\mathfrak{U}_N=\{u\in \mathfrak{U}: \; |u(t,\omega)|\leq N, \; \forall t\geq 0, \; \omega \in \Omega\},
$$
and the associated value function
$$
V_N(x)=\sup_{u\in \mathfrak{U}_N} \mathbb{E} \left[ \int \limits_0^{\infty} e^{-\rho t} \left( \ln u(t)-cx^2(t) \right)dt \right]
$$		
where $x(\cdot)$ is the solution to \eqref{sldyn} with $x(0)=x$.

\begin{lemma}\label{Napprox}
	For all $x\geq 0$ we have
	$$
	0\leq V(x)-V_N(x)\leq \frac{(\rho+b)^2}{4\rho cN^2}
	$$
\end{lemma}
\noindent
For $\sigma>0$, the ellipticity of the HJB equation \eqref{OHJB} in $(0, \infty)$ induces extra regularity for the function V in $(0, \infty)$. Hence, the optimal dynamics for the shallow lake problem are described by
\begin{equation}\label{opt}
\begin{cases}
	dx_*(t)=\left(-\frac{1}{V_\sigma'(x_*(t))}-bx(t)+r\big(x_*(t)\big)\right)dt+\sigma x_*(t)dW_t & \\
	x_*(0)=x & 
\end{cases}
\end{equation}
\noindent
Standard results in diffusion theory (e.g., Lemma 23.18 in \cite{kallenberg}) imply that the process $\{x_*(t):\ t\ge 0\}$ has an invariant distribution with density $f$. In particular, for any initial state, $x\ge 0$, and any Borel subset of $\RR$, $A$, we have
\[
\lim_{t\to\infty} \PP\big[x_*(t)\in A\big]=\int_A f(x)\ dx.
\]
The last proposition describes the invariant distribution of $x_*$ and the precise asymptotics of its tails, as $x\to 0$ and as $x\to\infty$.
\noindent
\begin{propos}\label{asympt}
	The density, $f$, of the stationary distribution of the optimal dynamics \eqref{opt} is
		\begin{equation}\label{inv1}
		f(x)=\frac{1}{Z} x^{-2\big(1+\frac{b}{\sigma^2}\big)}\, e^{-\frac{2}{\sigma^2}\Phi_\sigma (x)},
	\end{equation}
	where $Z$ is a normalising constant and 
	\[
	\Phi_{\sigma}(x)=\int \limits_{x}^{\infty}\Big(-\frac{1}{V_\sigma'(u)} + {r(u)}\Big)\frac{du}{u^2}, \qquad x>0.
	\]
	In particular,	\begin{equation}\label{tails}
		\lim_{x\to 0} x\ \Phi_{\sigma}(x) = \frac{1}{ |V_\sigma'(0)|}, \qquad\text{and}\qquad  \lim \limits_{x\rightarrow \infty}\Phi_\sigma (x)=0.
	\end{equation}
\end{propos}
	
		\subsection{Proofs}\label{proofs}
		
In this subsection, we present the proofs that are not straightforward modifications of those in  \cite{KLS}. 	\\
	
\noindent
{\bf Proof of Proposition \ref{V0fin}}\\ Based on Proposition \ref{propdyn}(i), for any $u\in\mathfrak{U}$,  $x(t)\ge M_t(u)$. In addition, for any positive $\mathbb{P}$-a.s. locally integrable $\mathcal{F}_t$-adapted process $f$, by Lemma A.1(i) in \cite{KLS}) we have that
\begin{equation}\label{Mt}
\mathbb{E}\left[\int \limits_0^{\infty} e^{-\rho t}M_t(f)dt\right]=\frac{1}{\rho+b}\mathbb{E}\left[\int \limits_0^{\infty}e^{-\rho t}f(t)dt\right].
\end{equation}
\noindent
Therefore, using first Jensen's  inequality and then \eqref{Mt}, we find
\begin{align*}
	\frac{\ln c}{2\rho}+\mathbb{E}\left[\int_0^\infty e^{-\rho t}\ln u(t)\ dt\right] &\le\frac{1}{\rho}\ln\mathbb{E}\left[\int_0^\infty \rho e^{-\rho t}\sqrt{c}u(t)\ dt\right] \\
	&=\frac{1}{\rho}\ln\mathbb{E}\left[\int_0^\infty \rho(\rho+b) e^{-\rho t}\sqrt{c}M_t(u)\ dt\right] \\
	&\le \frac{\ln(b+\rho)}{\rho}+\frac{1}{\rho}\ln\mathbb{E}\left[\int_0^\infty \rho e^{-\rho t}\sqrt{c}x(t)\ dt\right] \\
	&\le \frac{\ln(b+\rho)}{\rho}+\frac{1}{2\rho}\ln\mathbb{E}\left[\int_0^\infty \rho e^{-\rho t}cx^2(t)\ dt\right].
\end{align*}

In view of (\ref{ac_constraint}), it suffices to consider  $u\in\mathfrak{U}$ such that $D=\mathbb{E}\left[\int_0^\infty e^{-\rho t}cx^2(t)\ dt\right]<\infty$.
Then,
\[
\mathbb{E} \left\{ \int_{0}^{\infty}e^{-\rho t}\big[\ln u(t)-cx^{2}(t)\big]dt \right\}\le  \frac{\ln(b+\rho)}{\rho}+\frac{\ln(\rho D)}{2\rho}-D-\frac{\ln c}{2\rho}\le\frac{1}{\rho}\ln\left(\frac{b+\rho}{\sqrt{2ec}}\right).
\]
The assertion follows by taking the supremum over admissible controls.\hfill$\Box$	\\

\noindent
{\bf Proof of Proposition \ref{vprop1}} \\
\noindent
{\em Proof of the  lower bound:}\quad 
It suffices to produce an admissible control that achieves a benefit $J(x;u)$ greater than or equal to the bound. For $x,y\in\RR$, let us denote by $x\vee y=\max\{x,y\}$ and consider the feedback control 
\begin{equation}\label{feedcon}
u(t)=\frac{1\vee x(t)}{1+x^2(t)}+a-r(x(t)).
\end{equation}
An elementary variation of parameters argument for stochastic differential equations yields that
\begin{align}\label{xZ}
x(t)&=xZ_t+\int \limits_0^t \frac{Z_t}{Z_s} \left(u(s)+r(x(s)) \right)ds \\ &=xZ_t+\int_0^t \frac{Z_t}{Z_s}\ \Big(a+\frac{1\vee x(t)}{1+x^2(t)}\Big)ds\nonumber \\
&=xZ_t+aM_t(1)+\int_0^t \frac{Z_t}{Z_s}\ \frac{1\vee x(s)}{1+x^2(s)}ds.
	\label{testX}
\end{align}
We now square both sides to get
\begin{align*}
	x^2(t)&=x^2Z_t^2+2axZ_tM_t(1)+a^2M_t^2(1)+\left(\int_0^t \frac{Z_t}{Z_s}\ \frac{1\vee x(s)}{1+x^2(s)}ds\right)^2\\
	&\quad+2aM_t(1)\int_0^t \frac{Z_t}{Z_s}\ \frac{1\vee x(s)}{1+x^2(s)}ds+2\int_0^t \frac{Z_t^2}{Z_s}\ \frac{x \cdot 1\vee x(s)}{1+x^2(s)}ds.
\end{align*}
For the fourth and fifth term on the right-hand-side we simply use the inequality $\frac{1\vee x(s)}{1+x^2(s)}\le 1$. For the rightmost term we note that by (\ref{testX}) we have  $x\le x(s)Z_s^{-1}$. It follows that
\begin{equation}
	x^2(t)\le x^2Z_t^2+2axZ_tM_t(1)+(a+1)^2M_t^2(1)+2\int_0^t \frac{Z_t^2}{Z_s^2}\ ds.
	\label{testX2}
\end{equation}
It is straightforward to check that $\mathbb{E}\big[Z_t^2\big]=\exp((\sigma^2-2b)t)$, while by Lemma A.1 (ii) in \cite{KLS}, we have
\begin{equation}\label{ZMint}
\mathbb{E}\big[\int_0^\infty e^{-\rho t}Z_tM_t(1)dt\big]=\frac{A}{c(\rho+b)}.
\end{equation}
Hence, we can find  some constant $B$, such that
\begin{equation}
	\int_0^\infty e^{-\rho t}\EE\big[cx^2(t)\big]\ dt \le A\big(x+\frac{a}{\rho+b}\big)^2+B.
	\label{Xener}
\end{equation}

On the other hand, using that for all $x\ge 0$, we have $a\ge r(x)$ and $\displaystyle{\frac{1\vee x}{1+x^2}\ge \frac{1}{1+x} }$,  we find
\begin{align*}
	\int_0^\infty e^{-\rho t}\EE\big[\ln u(t)\big]\ dt &\ge-\int_0^\infty e^{-\rho t}\EE\big[\ln\big(1+x(t)\big)\big]\ dt \\
	&\ge -\frac{1}{\rho}\ln\left(\int_0^\infty \rho e^{-\rho t} \Big(1+\EE\big[x(t)\big]\Big) dt\right)\\
	&=-\frac{1}{\rho}\ln\left(1+\rho\int_0^\infty  e^{-\rho t} \EE\big[x(t)\big] dt\right),
\end{align*}
where in the penultimate step we have used Jensen's inequality.\\ 
By  (\ref{testX}) it follows that
\[\displaystyle
\EE\big[x(t)\big]\le x\EE\big[Z_t\big]+(a+1)\EE\big[M_t(1)\big]=xe^{-bt}+(a+1)\EE\big[M_t(1)\big].
\]
Hence, using Lemma A.1(i) in \cite{KLS}, we obtain that the control $u$ in \eqref{feedcon} satisfies
\[
\int_0^\infty e^{-\rho t}\EE\big[\ln u(t)\big]\ dt \ge -\frac{1}{\rho}\ln\left(1+\frac{\rho x}{\rho+b}+\frac{a+1}{\rho+b}\right).
\]

The preceding estimate and (\ref{Xener}) together imply that, for some suitable constant $K_1$,
\begin{align*}
	V(x)&\ge J(x;u) =\EE\big[\int_0^\infty e^{-\rho t}\big(\ln u(t)-cx^2(t)\big)\ dt \big]\\
	&\ge -A\left(x+\frac{a}{b+\rho}\right)^2-\frac{1}{\rho}\ln \left(x+\frac{a}{b+\rho}\right)+K_1.
\end{align*}

\vskip.07in	
\noindent
{\em Proof of the upper bound:}\quad Fix $u\in\mathfrak{U}$. Equation \eqref{xZ} gives
\begin{align*}
	x^2(t)&\ge x^2Z_t^2+2xZ_t^2\int_0^t\frac{1}{Z_s}\left(u(s)+r(x(s))\right) ds\\&=x^2Z_t^2+2xZ_tM_t(a+u)-\int_0^t\frac{Z_t^2}{Z_s}2x(a-r(x(s)) ds\\
	&\geq x^2Z_t^2+2axZ_tM_t(1)+2xZ_tM_t(u)-\int_0^t\frac{Z_t^2}{Z_s^2}2x(s)(a-r(x(s)) ds.
\end{align*}

Since
\[
\lim \limits_{x\rightarrow \infty} x(a-r(x))=C\in \mathbb{R} \; \Rightarrow \; x(a-r(x))\leq K\; \text{for some positive constant K,} 
\]

we can further estimate $x^2(t)$ from below by
\begin{align}\label{mtuest}
	x^2(t)&\ge x^2Z_t^2+2axZ_tM_t(1)+2xZ_tM_t(u)-2K\int_0^t \frac{Z_t^2}{Z_s^2} ds.
\end{align}

Using the elementary inequality that $\ln a\le ab-\ln b-1$, which holds for all $a,b>0$, and Lemma A.1(ii) in \cite{KLS}, we obtain, for some constant B,
\begin{align*}
	&\int_0^\infty e^{-\rho t}\EE\big[\ln u(t)\big]dt \le \EE\left[\int_0^\infty e^{-\rho t} \Big\{2AxZ_tu(t)-\ln\big(2AxZ_t\big)\Big\}\ dt\right]\nonumber\\
	&\qquad=\EE\left[\int_0^\infty\hspace{-2mm} e^{-\rho t} 2cxZ_tM_t(u)\, dt\right]-\frac{\ln(2Ax) }{\rho}+\frac{2b+\sigma^2}{2\rho^2}\nonumber\\
	&\qquad\le \EE\left[\int_0^\infty \hspace{-2mm}e^{-\rho t} \Big(cx^2(t)-cx^2Z_t^2-2acxZ_tM_t(1)+2Kc\int_0^t\frac{Z_t^2}{Z_s^2}\, ds\Big)\, dt\right]-\frac{\ln x + B}{\rho}\nonumber,
\end{align*}
where in the last step we have used \eqref{mtuest} to estimate $Z_tM_t(u)$. By another application of \eqref{ZMint}  we conclude that for every $u\in\mathfrak{U}$ there exists  $K_2>0$ such that
\[
J(x;u)\le -A\left(x+\frac{a}{b+\rho}\right)^2-\frac{1}{\rho}\ln \big(x+\frac{a}{\rho+b}\big)+K_2.
\]

The assertion now follows by taking the supremum over $u\in\mathfrak{U}$. \hfill$\Box$\\[2mm]

\noindent
{\bf Proof of  Proposition \ref{vprop2}} \\
By Assumption \ref{r}, there exists $\eps>0$ such that $r(x)<(b+\rho)x,\; \forall x\in\; (0,\eps]$. In view of Proposition \ref{Vdec}, it suffices to assume that $x_2\le \eps$, since otherwise we have
\[
V(x_2)-V(x_1)\le -A(x_2^2-x_1^2)< -A\eps(x_2-x_1).
\]

For a positive constant $d$, choose  a $u_d\in\mathfrak{U}$ that is constant and equal to $d$ up to time $\tau_d=\tau_{u_d}$.
Then, Proposition \ref{specialdpp} yields
\[
V(x_1)\ge \frac{\ln d-cx_2^2}{\rho}(1-\EE\big[e^{-\rho\td}\big])+\EE\big[{e^{-\rho\td}}\big] V(x_2),
\]
or equivalently,
\begin{equation}
	\big(V(x_2)-V(x_1)\big)\EE\big[e^{-\rho\td}\big]\le -\big(\ln d-\rho V(x_1)-cx_2^2)\ \EE\left[\int_0^{\td}e^{-\rho t}\, dt\right].
	\label{ub}
\end{equation}

Consider now the solution $x_d(\cdot)$ to (\ref{sldyn}) with $x(0)=x_1$ and control $u_d$. We can apply It\^o's formula to $e^{-\rho t}x_d(t)$, 
followed by the optional stopping theorem for the bounded stopping time $\tau_N=\td\wedge N$, and we get
\begin{align}
	\EE\big[e^{-\rho \tau_N}x_d(\tau_N)\big]-x_1&=\EE\left[\int_0^{\tau_N} e^{-\rho t}\big(d-(b+\rho)x_d(t)+r(x_d)\big)\, dt\right].
	\label{lastone}
\end{align}

The leftmost term of \eqref{lastone} is equal to \(x_2\EE\big[e^{-\rho \td};\td\le N\big]+e^{-\rho N}\EE\big[x_c(\tau_N); \td>N\big]\). 
\vskip.05in

On the other hand, since we have assumed that $x_2\le \eps$, we have $x_d(t)\le \eps$ up to time $\td$. Thus, the right hand side 
of \eqref{lastone} is bounded by \( \EE\left[\int_0^{\tau_N}e^{-\rho t} d\ dt\right]\). 
\vskip.05in

Letting $N\to\infty$ in $\eqref{lastone}$, by the monotone convergence theorem,  we have
\[
x_2 \EE\big[e^{-\rho \td}\big]-x_1\le d\,\EE\big[\int_0^{\td} e^{-\rho t}dt\big] \Longleftrightarrow    (x_2-x_1)\EE\big[e^{-\rho\tau_d}\big]\le (d+\rho x_1)\ \EE\big[\int_0^{\td} e^{-\rho t}dt\big].
\]

Substituting this in (\ref{ub}) and choosing $\ln d=\rho V(x_1)+1+cx_2^2$, we find
\begin{equation}\label{lowb}
	V(x_2)-V(x_1)\le -(x_2-x_1)\left(e^{\rho V(x_1)+1+cx_2^2}+\rho x_1\right)^{-1} \!\!\!.
\end{equation}

The assertion now follows setting $c(x_2)=A\eps \mathbf{1}\{x_2>\eps\}+ \left(e^{\rho V(0)+1+cx_2^2}+\rho x_2\right)^{-1} \mathbf{1}\{x_2\leq \eps\}  $
and 
$C_1=A\eps \wedge \left(e^{\rho V(0)+1+c\eps^2}+\rho \eps\right)^{-1}>0$. \hfill$\Box$ \\[2mm]
\vskip.07in	
		\noindent
{\bf Proof of Theorem \ref{asymptotic} }\\We define an auxiliary function $v$ by
\begin{equation*}
	v(x):=V(x)+A\left(x+\frac{a}{b+\rho}\right)^2+\frac{1}{\rho}\ln \left(2A(x+\frac{a}{b+\rho})\right)-K, \quad x\in\RR. 
\end{equation*}

Straightforward calculations yield that   $v$ is a viscosity solution in $(0, \infty)$ of the equation
\begin{multline}
	\rho v+\left(bx-r(x)\right)v'+\ln\left(1+\frac{1-\rho\big(x+\frac{a}{b+\rho}\big)v'}{2A\rho\left(x+\frac{a}{b+\rho}\right)^2} \right)-\dfrac{1}{2}\sigma^{2}x^{2}v''+f=0,
\end{multline}
where 
\[
f(x)=\frac{a(b+\frac{\sigma^2}{2})+(b+\rho)r(x)}{\rho\big(a+x(b+\rho)\big)}+\frac{\sigma^2x(b+\rho)}{2\rho\big(a+x(b+\rho)\big)^2}-\frac{2A}{b+\rho}(a-r(x))(a+x(b+\rho))+2AC.
\]
Note that $f$ is smooth on $[0,\infty)$ and vanishes as $x\to\infty$.
\vskip.05in

Let $v_\lambda(y)=v(\frac{y}{\lambda})$ and observe that, if $v_\lambda(1) \rightarrow 0$ as $\lambda \rightarrow 0$, then $v(x) \rightarrow 0$ as $x 
\rightarrow \infty$. It turns out that  $v_\lambda$ solves
\begin{multline*}
	\rho v_\la+ \left(bx-\lambda r(\frac{x}{\lambda})\right)v'_\la+\ln\left(1+\frac{\la^2\big(1-\rho\big(x+\frac{\la a}{b+\rho}\big)v'_\la\big)}{2A\rho\left(x+\frac{\la a}{b+\rho}\right)^2} \right)-\dfrac{1}{2}\sigma^{2}x^{2}v''_\la+f\big(\frac{x}{\la}\big)=0.
\end{multline*}

\noindent
Since, by (\ref{Vbounds}) $v_\la $ is uniformly bounded, we consider the half-relaxed limits $v^*(y)=\limsup_{x \rightarrow y, \la\rightarrow 0}v_\la(x)$ and $v_*(y)=\liminf_{x \rightarrow y, \la\rightarrow 0}v_\la(x)$ 
in $(0, \infty)$. By \cite{BP}, $v^*$ and $v_*$ are respectively sub- and super-solutions of 
\begin{equation}
	\rho w +byw'-\frac{1}{2} \sigma^2 y^2 w''=0.
\end{equation}

It is easy to check that for any $y>0$ we have $v^*(y)=\limsup_{x\to\infty}v(x)$ and $v_*(y)=\liminf_{x\to\infty}v(x)$. 
\vskip.05in
The subsolution property of $v^*$ and the supersolution property of $v_*$ give 
\[
\limsup_{x\to\infty}v(x)\le 0\le\liminf_{x\to\infty}v(x)\le \limsup_{x\to\infty}v(x).
\]\qed

\noindent
{\bf Proof of Lemma \ref{Napprox}: }
	The lower bound is evident since the supremum in the definition of $V_N$ is taken over a subset of admissible controls. Let now $\varepsilon >0$, and consider a control $u\in \mathfrak{U}$ such that
	$$
	V(x)<\mathbb{E}_x \left[ \int \limits_0^{\infty} e^{-\rho t} (\ln u(t) -c x^2(t))dt \right]+\varepsilon
	$$
	Let us now define $u_N=u\wedge N=\min\{u,N\}\in \mathfrak{U}_N,$ and $x_N(\cdot)$ the solution to \eqref{sldyn} with $u_N$ as control. Clearly,
	$$
	V_N(x)\geq \mathbb{E}_x \left[ \int \limits_0^{\infty}e^{-\rho t} (\ln u_N(t) -c x_N^2(t))dt  \right]
	$$
	so that
	\begin{equation}\label{VVN}
		V(x)-V_N(x)\leq \varepsilon +\mathbb{E}_x \left[ \int \limits_0^{\infty}e^{-\rho t} (\ln\left( \frac{u(t)}{u_N(t)}\right) -cx^2(t)+cx_N^2(t))dt  \right]. 
	\end{equation}
	If we denote $\Delta u=(u-N)^+=\max\{u-N,0\}$, we can write $u=u_N+\Delta u$,
	with $\Delta u \neq 0$ if and only if $u>N$. Thus,
	
	\begin{multline}\label{lnDu}
		\mathbb{E} \left[ \int \limits_0^{\infty}e^{-\rho t} \ln\left( \frac{u(t)}{u_N(t)}\right)dt \right]=\mathbb{E} \left[ \int \limits_0^{\infty}e^{-\rho t} \ln\left(1+ \frac{\Delta u(t)}{N}\right)dt \right] \\ \leq \frac{1}{N} \mathbb{E} \left[ \int \limits_0^{\infty} e^{-\rho t} \Delta u(t)dt \right]=\frac{\rho+b}{N} \mathbb{E} \left[ \int \limits_0^{\infty} e^{-\rho t} M_t( \Delta u )dt \right],
	\end{multline}
	\noindent
	where the final equality is due to \eqref{Mt}. After application of \eqref{xZ} for both controls $u$ and $u_N$, we get
	$$
	x(t)-x_N(t)=M_t(\Delta u )+\int \limits_0^t \frac{Z_t}{Z_s} \left( r(x(s))-r(x_N(s)) \right)ds. 
	$$
	 In view of Propositions \ref{positive} and \ref{monotone}, we have $0\leq x_N(\cdot) \leq x(\cdot),$ $\mathbb{P}$-a.s. and since $r$ is increasing, we have
	$$
	x^2(t)-x^2_N(t)\geq (x(t)-x_N(t))^2\geq M^2_t(\Delta u), \; \forall t \geq 0, \mathbb{P}-\text{a.s.}
	$$
	The preceding estimate and \eqref{lnDu} can further strengthen inequality \eqref{VVN} to
	$$
	V(x)-V_N(x)\leq \varepsilon + \mathbb{E} \left[ \int \limits_0^{\infty} e^{-\rho t}\left( \frac{\rho+b}{N} M_t(\Delta u)-c M^2_t(\Delta u) \right)dt  \right]\leq \varepsilon +\frac{1}{\rho} \left(\frac{\rho+b}{2N\sqrt{c}}\right)^2.
	$$
	Since $\varepsilon>0$ was arbitrary, the assertion follows.\qed\\
	
	\noindent
	{\bf Proof of Lemma \ref{asympt} }\\
	To calculate the stationary density $f$ of an It\^{o} diffusion $y_t$, we may solve the stationary Fokker-Planck equation,
	\begin{equation}\label{FPe}
		\mathcal{L}^*(f)=0, 
	\end{equation}
	where $\mathcal{L}^*$ is the adjoint of the generator $\mathcal{L}$ of the process, and demand that the solution satisfy the constraints $f\geq 0$ and  $\int f=1$.
	
	Equation  \eqref{FPe} takes a convenient form, when the corresponding process $y_t$ has a constant diffusion coefficient. Specifically, if $dy_t=-g(y_t)dt+\sigma dW_t$, then $\mathcal{L}(f)=\frac{\sigma^2}{2}f''- gf'$ and equation \eqref{FPe} becomes 
	\begin{equation}
		\frac{\sigma^2}{2} \frac{d^2f}{dy^2}+\frac{d(gf)}{dy}=0.
	\end{equation} 
	Hence, if $\{y_t: t\ge 0\}$ has a unique invariant distribution, its density, $f_Y$, will be given by
	\begin{equation}\label{invdist}
		f_Y(y)=\frac{1}{Z} e^{-\frac{2}{\sigma^2}G(y)},\quad y\in\RR,
	\end{equation} 
	where $G$ is an antiderivative of $g$, and $Z$ is a normalising constant.
	We can reduce the dynamics of the optimally controlled stochastic lake to the preceding form via the transformation $y_t=\ln (x_{*}(t))$. Indeed, applying It$\hat{o}$'s rule to the process $y_t=\log x_{*}(t)$, we find 
	\begin{equation}\label{con_dif}
		dy_t=\left(e^{-y_t}h(e^{y_t})-\frac{\sigma^2}{2}\right)dt+\sigma dW_t
	\end{equation}
	where $h(x)=-\frac{1}{V'(x)}-bx+r(x)$ is the drift of the optimally controlled lake, $x_*$. Assumption \ref{r} and Proposition \ref{vprop}(ii) ensure that $y$ in  \eqref{con_dif} has a unique invariant distribution, whose density is given by \eqref{invdist}, with
	\[
	G(y)=\big(b+\frac{\sigma^2}{2}\big)y + \Phi_\sigma(e^y).
	\]
	Hence, $x_*(t)=e^{y_t}$ itself has a unique invariant distribution, whose density is given by
	\begin{equation}
	f(x)=\frac{f_Y(\ln x)}{x}=\frac{1}{Z} x^{-2\big(1+\frac{b}{\sigma^2}\big)}\, e^{-\frac{2}{\sigma^2}\Phi_\sigma (x)},\quad x>0.
	\end{equation}
	The limits in \eqref{tails} that determine the tail asymptotics of the invariant distribution are an immediate consequence of Assumption \ref{r} and Proposition \ref{vprop}(ii).\qed\\
	
		\section{Numerical Investigation}\label{numerics}
		In this section we implement the monotone numerical scheme suggested in \cite{KLS} to numerically compute the value function of \eqref{ihvf}. We do so for parameter values in ranges that correspond to distinct qualitative behaviour, as well as for recycling rates approximating a step function. Computation of the value function provides access to the invariant distribution of the optimally controlled stochastic shallow lake \eqref{opt} through equation \eqref{inv1} and we investigate the shape of the invariant density. Finally, we can simulate paths of \eqref{opt} and explore the statistics of the transition times between oligotrophic and eutrophic states.	
		
		More precisely, we consider a computational domain $[0,l]$, for a sufficiently large $l$, and a uniform partition $0 = x_0 < x_1 < \ldots < x_{N-1} < x_N = l$, whose size we denote by $\Delta x$. Having in mind (\ref{DVbounds}),  if $V_{i}$ is the approximation of $V$ at $x_i$,  we employ a backward finite difference discretisation  to approximate the first derivative in the linear term of the (OHJB), a forward finite difference discretisation for the derivative in the logarithmic term and a central finite difference scheme to approximate the second derivative. These considerations yield, for $i = 1,\ldots, N-1$,   the approximate equations
		\begin{multline}\label{DHJB}
			V_i - \frac{1}{\rho} \Big(r(x_i) - bx_i \Big) \frac{ V_{i} - V_{i-1} }{ \Delta x} + \frac{1}{\rho}
			\left[cx_i^2 + 1 +\ln \left( - \frac{V_{i+1} - V_i}{\Delta x}\right)\right] \\
			- \frac{\sigma^2}{2\rho}x_i^2
			\frac{ V_{i+1} + V_{i-1} - 2 V_i}{(\Delta x)^2}=0.
		\end{multline}
		
		In addition, we impose a boundary condition at the right endpoint, $l$, which we assumed to be sufficiently large, taking advantage of the asymptotic behaviour of $V$  described in Theorem \ref{asymptotic}. That is, 
		\begin{equation}\label{rboundary}
			V_N=-A\left(l+\frac{a}{b+\rho}\right)^2-\frac{1}{\rho}\ln \left(2A(l+\frac{a}{b+\rho})\right)+K.
		\end{equation}
		
		Finally, another approximate equation is provided by the boundary condition at $x=0$, as given in Corollary \ref{coro1}, i.e.,
		\begin{equation}\label{lboundary}
		V_0+\frac{1}{\rho}
				\left[1 +\ln \left( - \frac{V_{1} - V_0}{\Delta x}\right)\right]=0.   
		\end{equation}
		
		In fact, equation \eqref{lboundary} coincides with \eqref{DHJB} for $i=0$, and we have in total a system of $N$ nonlinear equations for the $N$ unknowns, $V_0,\ldots,V_{N-1}$. In this paper, we approximate the solution of this system using the Newton-Raphson method. As an initial estimation of the solution in the Newton-Raphson algorithm, we considered a quadratic function $V^0$, such that $V^0(x_N)=V_N,\; V^0(x_0)=V^0(0)=\frac{1}{\rho}\ln\left(\frac{b+\rho}{\sqrt{2ec}}\right) $ (the upper bound in Proposition \ref{V0fin}) and $V^0(x_1)=V^0(0)-\Delta x e^{-(\rho V^0(0)+1)}$, so that \eqref{lboundary} is initially satisfied.

		The proof of convergence of the numerical scheme to the value function follows a general argument proposed in \cite{Barles} to prove convergence of monotone schemes for viscosity solutions of fully nonlinear second-order elliptic or parabolic, possibly degenerate, partial differential equations. By Proposition 2 in \cite{KLS}, the numerical scheme described above is consistent and monotone, provided
		\begin{equation}\label{num_condition}
				\Delta x \big(r(x) - bx\big)\leq \frac{\sigma^2}{2}.
			\end{equation}
In addition, as $\Delta x\to 0$, the continuous function that interpolates the solution at the points $x_0,\ldots,x_N$ converges locally uniformly to the unique constrained viscosity solution of \eqref{OHJB}, i.e., the actual value function. Condition \eqref{num_condition} is always satisfied for the stochastic shallow lake problem, provided $\Delta x$ is chosen suitably small. Note also that for the usual choice of recycling function, $r(x)=x^2/(x^2+1),$ condition \ref{num_condition} is satisfied if e.g. $b\geq 0.5$, even if $\sigma=0$.
	
			An advantage of this methodology for the computation of the value function is that we are free to choose any value of the parameter $\sigma$, as long as $\sigma^2<\rho+2b$ (see Proposition \ref{restr1}) and condition \ref{num_condition} are satisfied. In this way, we are not restricted to small values of the noise parameter $\sigma$. For instance, in Figure \ref{steady_iv} we extend to higher noise intensities the bifurcation diagram in Figure 5 of \cite{GKW}, revealing new features.
			
			In the first part of our numerical investigation, we study the problem with the typical choice of the function $r,$ i.e. $r(x)=x^2/(x^2+1)$, while in the second part we study the properties of the value function $V,$ which corresponds to a hyperbolic tangent function. 
			
			\subsection{The value function V and the optimal policy}
			
			In order to gain some first insight into the problem, we begin by exploring the properties of the value function $V$ for various values of the parameters $b,c,\rho, \sigma$. In our analysis, our choice of parameters is based on the bifurcation analysis in \cite{kslv}.
			
\begin{figure}
				\centering 
				\begin{subfigure}{0.5\textwidth}	
					\includegraphics[width=1\textwidth]{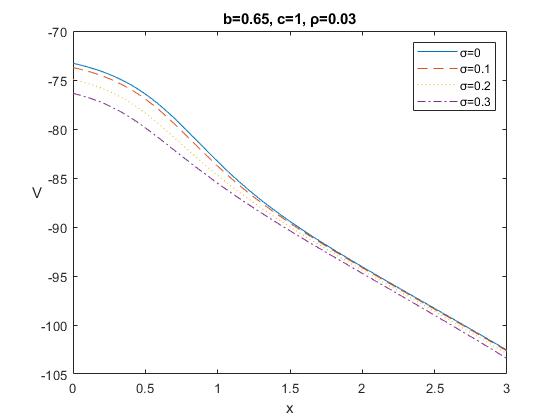}
					\caption{}
					\label{vfsec1}
				\end{subfigure}\hfil
				\begin{subfigure}{0.5\textwidth}
					\includegraphics[width=1\textwidth]{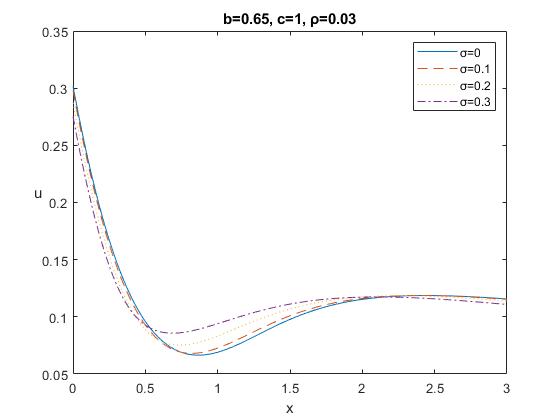}
					\caption{}
					\label{usec1}
				\end{subfigure}
				\begin{subfigure}{0.5\textwidth}
					\includegraphics[width=1\textwidth]{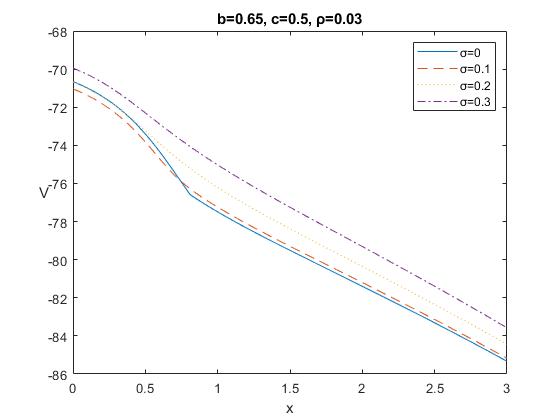}
					\caption{}
					\label{vf2sec1}
				\end{subfigure}\hfil
				\begin{subfigure}{0.5\textwidth}
					\includegraphics[width=1\textwidth]{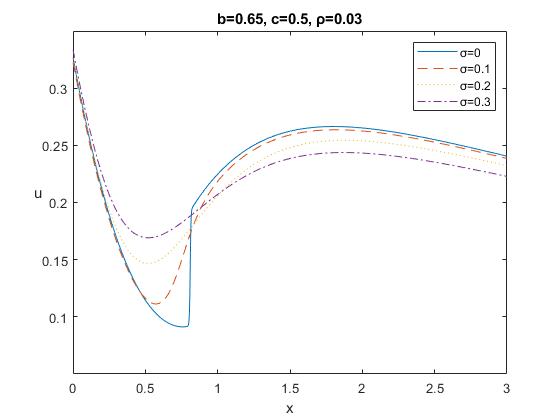}
					\caption{}
					\label{u2sec1}
				\end{subfigure}
				\vspace{-1\baselineskip}
				\caption{The value function $V$ (left) and the optimal policy (right) for different values of noise including the deterministic case $(\sigma=0)$. Parameters: Fig \ref{vfsec1}- \ref{usec1}$: (b,c,\rho)=(0.65,1,0.03)$ and Fig \ref{vf2sec1}-\ref{u2sec1}: $(b,c,\rho)=(0.65,0.5,0.03)$}.	
				\label{vf}
			\end{figure}

			Figures \ref{vfsec1} and \ref{vf2sec1} show the graph of the value function for the fixed parameters $(b,c,\rho)=(0.65,1,0.03)$ and $(b,c,\rho)=(0.65,0.5,0.03)$ respectively with the noise $\sigma$ varying. Notice that these graphs also depict the value function in the deterministic case ($\sigma =0$). In Figures \ref{usec1} and \ref{u2sec1}, the corresponding optimal management policies are shown. For the choice of parameters $(b,c,\rho)=(0.65,1,0.03)$, the optimal policies are smooth functions. 	On the other hand, when $(b,c,\rho)=(0.65,0.5,0.03)$,  the system exhibits a Skiba point and the (deterministic) optimal policy is discontinuous at this point. 
			
			\subsection{Invariant distribution}
			
			In this section, we numerically investigate the properties of the equilibrium distribution of the optimally controlled lake for different combinations of the parameters of the problem.

			Apart from the invariant density, $f,$ and cumulative distribution, $F,$ of the optimally controlled lake, we also present bifurcation diagrams based on its transformation invariant function, $I:=\sigma x f$. The main advantage of this object is its invariance under diffeomorphic coordinate transformations, which makes the transformation invariant function a suitable basis of bifurcation theory (see e.g. \cite{Zeeman} and \cite{GKW}). Following the definitions introduced in \cite{GKW}, the local maximisers of the transformation invariant $I$ are called \textit{stochastic attractors} of the process, while the local minimiser of $I$ is called the \textit{regime switching threshold}. The stochastic attractors are the natural analogue of the attracting steady states of the deterministic problem and the \textit{regime switching threshold} is the analogue of the indifference point (the Skiba point). 
			
\begin{figure}[H]
				\centering 
				\begin{subfigure}{0.5\textwidth}
					\includegraphics[width=1\textwidth]{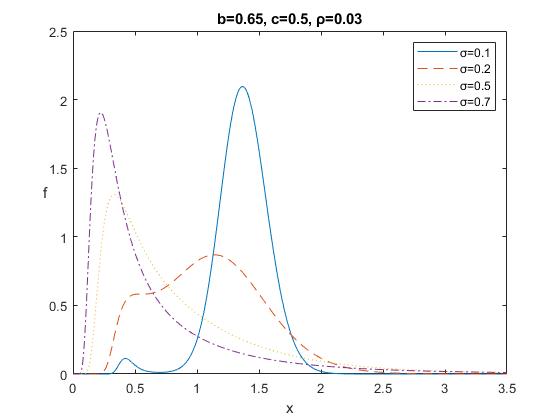}
					\caption{}
					\label{Iiv}
				\end{subfigure}\hfil
				\begin{subfigure}{0.5\textwidth}
					\includegraphics[width=1\textwidth]{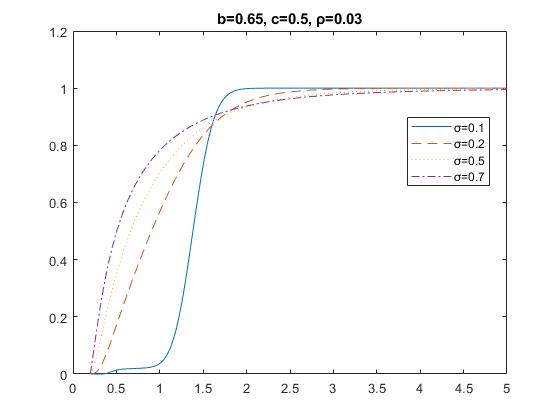}
					\caption{}
					\label{Fiv}
				\end{subfigure}
				\vspace{-1\baselineskip}
				\caption{Invariant density and cumulative distribution function for different values of the noise parameter $\sigma$. The choice of parameters is $(b,c,\rho)=(0.65,0.5,0.03)$.}	
				\label{IFiv}
			\end{figure}

\begin{figure}[H]
				\centering 
				\begin{subfigure}{0.5\textwidth}
					\includegraphics[width=1\textwidth]{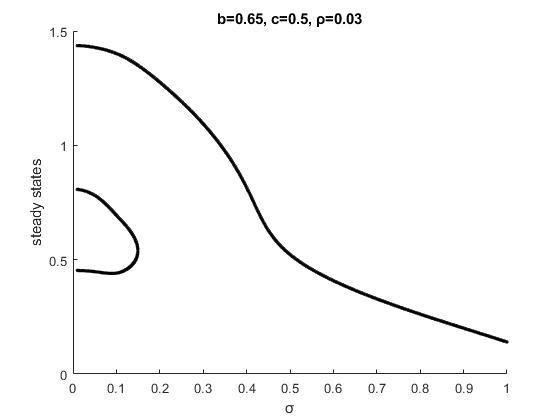}
					\caption{}
					\label{steady_iv}
				\end{subfigure}\hfil
				\begin{subfigure}{0.5\textwidth}
					\includegraphics[width=1\textwidth]{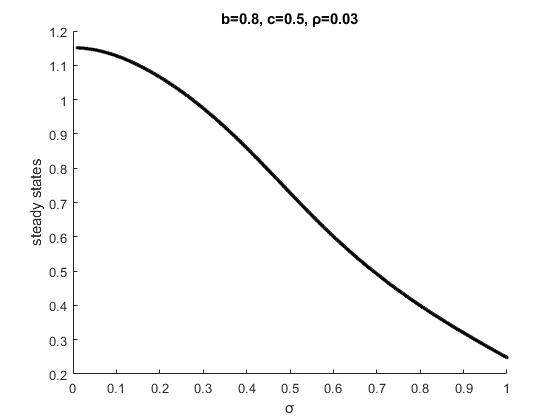}
					\caption{}
					\label{steady_viii}
				\end{subfigure}
				\vspace{-1\baselineskip}
				\caption{Bifurcation diagrams for the extrema of the transformation invariant function with respect to the noise parameter $\sigma$ when $(b,c,\rho)=(0.65,0.5,0.03)$ (left) and $(b,c,\rho)=(0.8,0.5,0.03)$ (right). The vertical axis corresponds to the location of the extrema of the transformation invariant function.}	
				\label{IFviii}	
			\end{figure}

			For the fixed parameters $(b,c,\rho)=(0.65, 0.5,0.03)$, Figure \ref{IFiv} shows the invariant density and cumulative distribution functions for several values of the noise parameter $\sigma$. For this set of parameters, the deterministic problem has a Skiba point. In the presence of small noise, the lake spends most of the time in the eutrophic state. As noise increases, the mode of the invariant distribution shifts to cleaner states, i.e. to lower concentrations of phosphorus. A detailed presentation of this shift from a bimodal distribution (with a peak at the eutrophic state) to a unimodal one (with a peak at oligotrophic phosphorus concentrations) due to the increase of noise is depicted in Figure \ref{steady_iv}. This bifurcation diagram illustrates the locations of the modes and the antimodes of the transformation invariant function, $I$, with respect to $\sigma$. 
			In the case of the fixed parameters $(b,c,\rho)=(0.8, 0.5,0.03)$, the deterministic problem exhibits a unique equilibrium in the eutrophic state (see \cite{W03}). Therefore, we have qualitatively different dynamics comparing to the preceding case.  In the presence of small noise, the (transformation) invariant function is unimodal with a peak at the eutrophic state, but the location of the mode moves to cleaner states as the noise intensity increases. These results are summarised in Figure \ref{steady_viii}. The same behaviour for large values of noise, as in the previous cases, is also present for combinations of parameters for which the deterministic problem exhibits a unique equilibrium in the oligotrophic state.  Based on the above observations, one could argue that noise seems to 'clean' the lake, in the sense that the lake spends more time in states corresponding to low phosphorus concentrations. One should bear in mind however that, at the same time, extremely polluted states become more likely at high noise intensities. Notice that if we were limited to small values of noise, e.g. $\sigma <0.2$ (see Figures \ref{steady_iv} and \ref{steady_viii}), we would not be able to observe this behaviour.
			
\begin{figure}[H]
				\centering 
				\begin{subfigure}{0.5\textwidth}
					\includegraphics[width=1\textwidth]{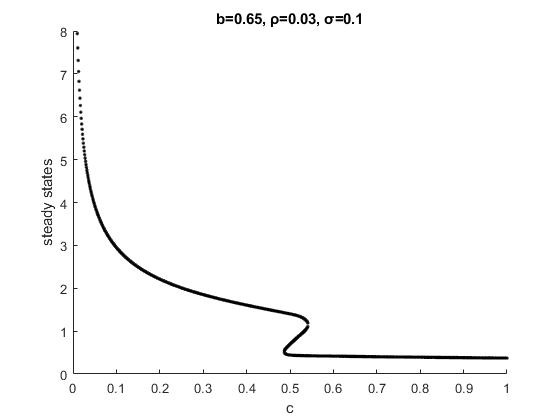}
					\caption{}
					\label{steady_c_065}
				\end{subfigure}\hfil
				\begin{subfigure}{0.5\textwidth}
					\includegraphics[width=1\textwidth]{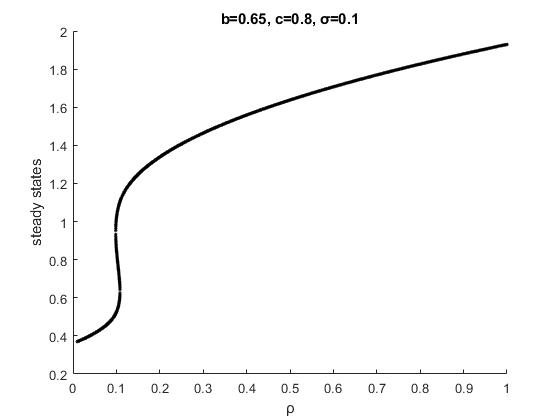}	
					\caption{}	
					\label{steady_r_065}	
				\end{subfigure}
				\vspace{-1\baselineskip}
				\caption{Left: Bifurcation diagram for the transformation invariant function with respect to the cost of pollution $c$ when $(b,\rho,\sigma)=(0.65,0.03,0.1)$. Right: Bifurcation diagram for the transformation invariant distribution with respect to the discount factor  $\rho$ when $(b,c,\sigma)=(0.65,0.8,0.1)$. The vertical axes correspond to the location of the extrema of the transformation invariant function.} 
				\label{steady_rc}
			\end{figure}

			Figure \ref{steady_c_065} illustrates a bifurcation diagram for the fixed parameters  $(b,\rho)=(0.65,0.03)$ and noise $\sigma=0.1$ with respect to the cost of pollution $c$. As it is expected by the definition of the total benefit \eqref{Jxu}, large values of $c$ attribute more weight to the ecological services of the lake, thus cleaning the optimally controlled lake. This is not the case, for the bifurcation diagram with respect to the discount factor $\rho.$ In Figure \ref{steady_r_065} the bifurcation diagram with respect to $\rho$ for the fixed parameters $(b,c)=(0.65,0.8)$ and noise $\sigma=0.1$ is depicted. In this diagram, we observe that as the discount factor $\rho$ increases and the benefit of future generations is discounted, the stochastic attractors of the system move towards eutrophic states.
			
			\subsection{The optimal paths and escape times} 
\begin{figure}
				\centering 
				\begin{subfigure}{0.5\textwidth}
					\includegraphics[width=1\textwidth]{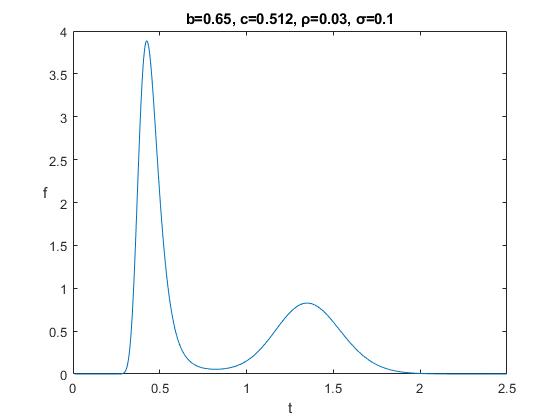}
					\caption{}
					\label{I0}
				\end{subfigure}\hfil
				\begin{subfigure}{0.5\textwidth}
					\includegraphics[width=1\textwidth]{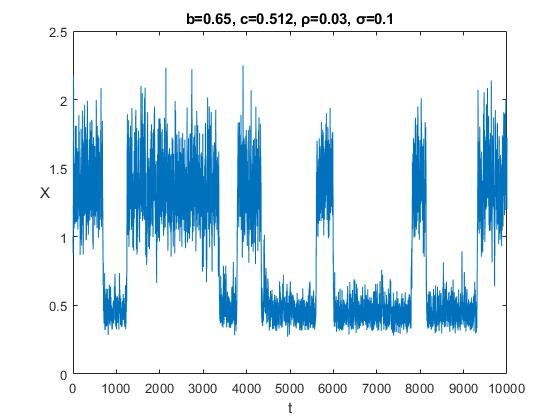}
					\caption{}
					\label{troxies201}
				\end{subfigure}
				\vspace{-1\baselineskip}
				\caption{The invariant density (left) and a simulated path (right) of the optimally controlled stochastic lake with two stochastic attractors. Parameters: $(b,c,\rho,\sigma)=(0.65,0.512,0.03,0.1)$}	
				\label{troxies2}
			\end{figure}

In the deterministic version of the problem, the optimally controlled system asymptotically approaches one of its attracting equilibrium states. On the other hand, noise introduces fluctuations around the stochastic attractors of the process. In case the stochastic system has more than one attractors, noise eventually induces fluctuations that are large enough to drive the system beyond the regime switching threshold and into the basin of attraction of a different attractor. The process thermalises there until a new large fluctuation causes another regime switching, and so on. The invariant density $f$ and one simulated path of the optimally controlled stochastic lake for the choice of parameters $(b,c,\rho, \sigma)=(0.65,0.512,0.03,0.1)$ are depicted in Figure \ref{troxies2}.
   		
			If we consider a diffusion $x_t$ in a double-well potential, $G$, with constant diffusion coefficient, i.e.,
			$dx_t=-G'(x_t)dt+\sigma dW_t$, and we denote by $x_{\pm}$ the stochastic attractors of the process, with $x_-<x_+$ and by $x_*$ the regime switching threshold, then the expected time $T_{x_+}$ of the system to hit $x_+$ when it starts at  $x_-$ is asymptotically exponential, in the sense that
			\begin{equation}\label{EXP}
				T_{x_+}/ \EE_{x_-} [T_{x_+}] \overset{d}{\rightarrow} Exp(1) \text{\; as \;}  \sigma \rightarrow 0
			\end{equation} 
			and it is described by the Arrhenius law
			\begin{equation}\label{AL}
				\displaystyle{\lim \limits_{\sigma \rightarrow 0} \frac{\sigma^2}{2}\log \EE_{x_-} [T_{x_+}]	=G_0(x_*)-G_0(x_-). }
			\end{equation}  
			The preceding results assume that the potential $G$ does not depend on the noise intensity, $\sigma$. This is not true in our case, as the optimal policy, hence the drift of the optimally controlled lake, depends on $\sigma$ through the HJB equation \eqref{OHJB}. The derivation of an expression analogous to \ref{AL} in this context is the subject of future work. \\
			
\begin{figure}[H]
				\centering
				\includegraphics[width=0.6\textwidth]{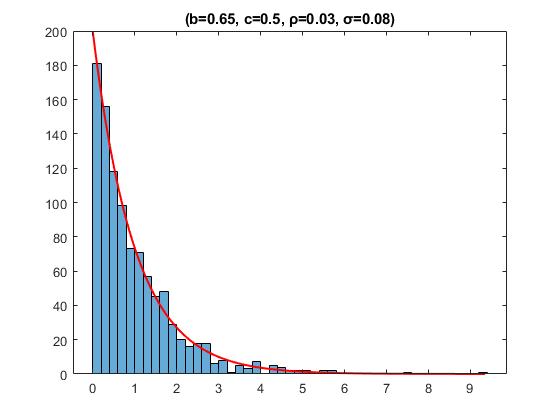}
				\caption{Histogram of 1.000 samples of the normalised transition time from the oligotrophic ($x_-$) to the eutrophic ($x_+$) state. Parameters: $(b,c,\rho,\sigma)=(0.65,0.5,0.03,0.08).$}	
				\label{esctime}
			\end{figure}

			Figure \ref{esctime} illustrates a histogram of 1000 realisations of the random variable $T_{x_+}/ \EE_{x_-} [T_{x_+}]$ for the optimally controlled lake. The expected time in the denominator was estimated by the sample mean of the computed times. The (red) curve corresponds to the exponential distribution with mean 1. The fit to the exponential distribution is very good for $\sigma = 0.08$.
			
			\subsection{The rate of recycling} 
\begin{figure}[H]
				\centering 
				\begin{subfigure}{0.5\textwidth}
			\includegraphics[width=1\textwidth]{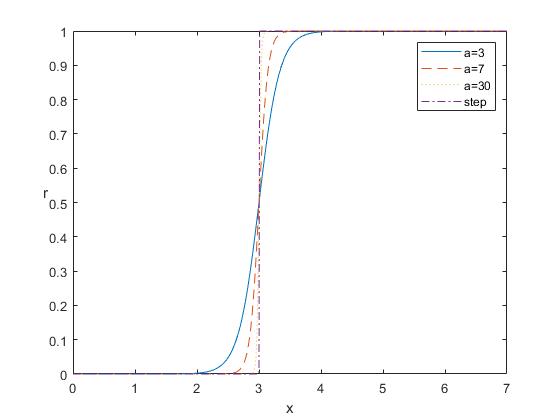}
				\vspace{-1\baselineskip}
				\caption{}	
				\label{tanh3}
				\end{subfigure}\hfil
				\begin{subfigure}{0.5\textwidth}
					\includegraphics[width=1\textwidth]{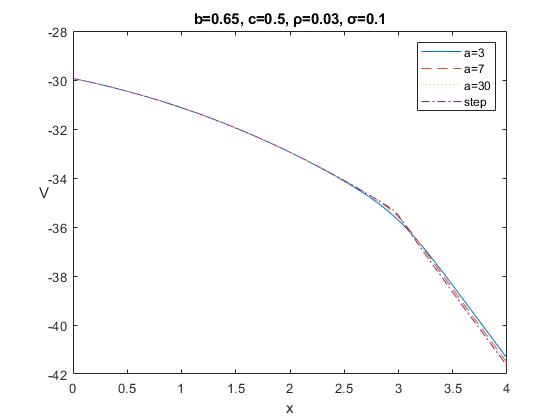}
					\caption{}
					\label{Va21}
				\end{subfigure}
				\vspace{-1\baselineskip}
				\caption{Left: The graph of the function $\frac{1}{2}(\tanh(a(x-3))+\tanh(3a))$ for different values of the parameter $a$ and the step function $\mathbf{1}\{x>3\}$. Right: the corresponding value function $V$ when $(b,c,\rho,\sigma)=(0.65,0.5,0.03,0.1)$}	
				\label{change}
			\end{figure}

\begin{figure}[H]
				\centering 
				\begin{subfigure}{0.5\textwidth}
					\includegraphics[width=1\textwidth,height=5cm]{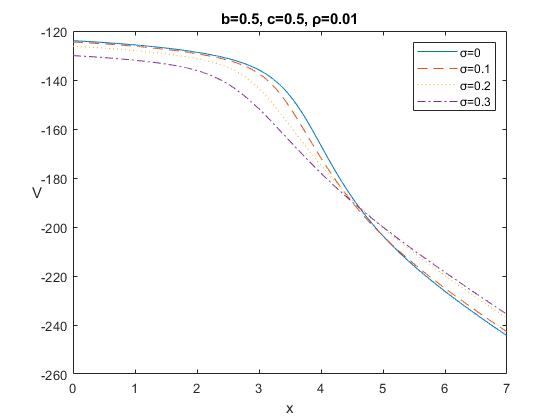}
					\caption{}
					\label{vfsec22}
				\end{subfigure}\hfil
				\begin{subfigure}{0.5\textwidth}
					\includegraphics[width=1\textwidth,height=5cm]{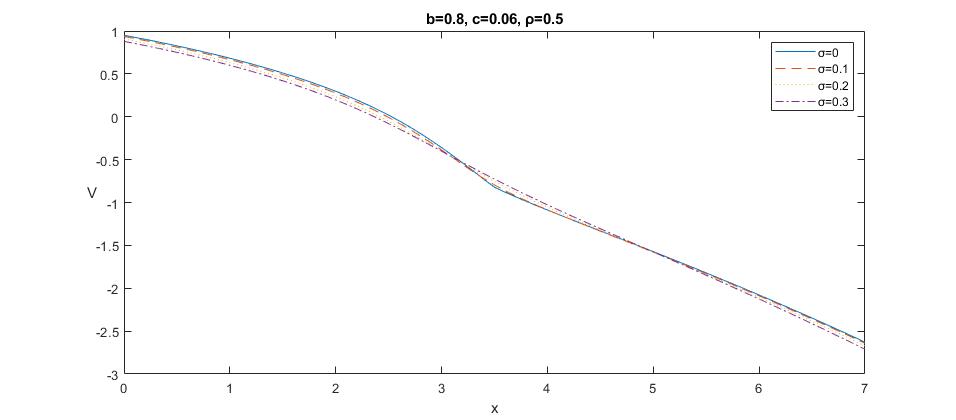}
					\caption{}
					\label{vfsec23}
				\end{subfigure}
				\begin{subfigure}{0.5\textwidth}
					\includegraphics[width=1\textwidth]
     {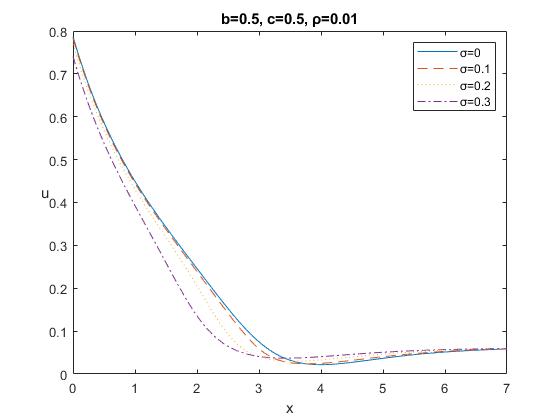}
					\caption{}
					\label{u22}
				\end{subfigure}\hfil
				\begin{subfigure}{0.5\textwidth}
					\includegraphics[width=1\textwidth]
     {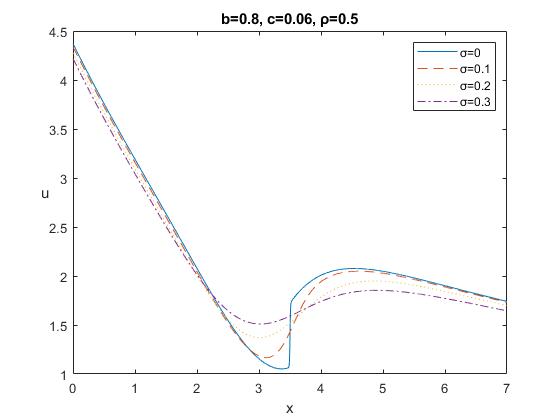}
					\caption{}
					\label{u23}
				\end{subfigure}
				\begin{subfigure}{0.5\textwidth}
					\includegraphics[width=1\textwidth]
     {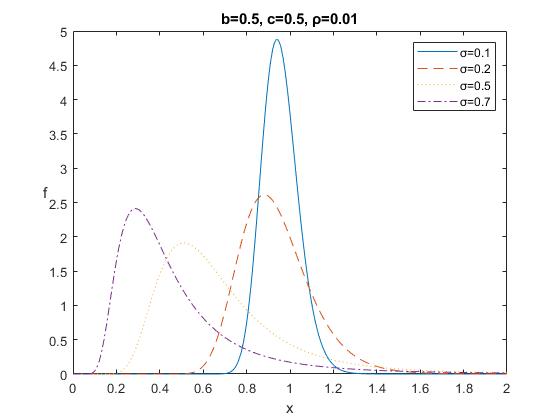}
					\caption{}
					\label{I22}
				\end{subfigure}\hfil
				\begin{subfigure}{0.5\textwidth}
					\includegraphics[width=1\textwidth]
     {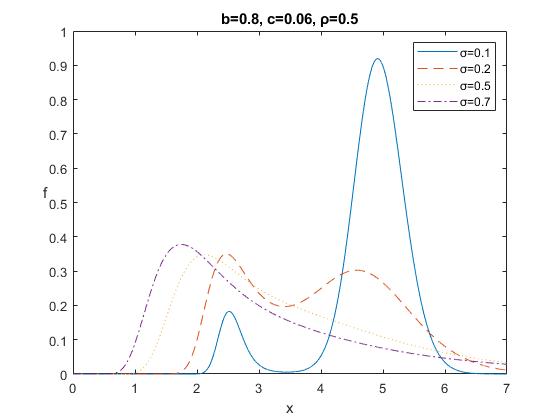}
					\caption{}
					\label{I23}
				\end{subfigure}
				\vspace{-1\baselineskip}
				\caption{The value function $V$, the optimal policy, $u$, and the equilibrium distribution, $f$ for different values of noise, when $(b,c,\rho)=(0.5,0.5,0.01)$ (left) and $(b,c,\rho)=(0.8,0.06,0.5)$ (right).}	
				\label{Iu21}
			\end{figure}

In this section, we present some numerical results, when a hyperbolic tangent function is used as the rate of recycling, $r$. We initially consider $r(x)=\tanh(x-3)+\tanh(3)$. In Figure \ref{Iu21}, we show the value function, the invariant density functions and the corresponding optimal policies for different combinations of the parameters $(b,c,\rho)$ and different values of $\sigma$. We observe that the lake has two attractors for small values of noise, when $(b,c,\rho)=(0.8,0.06,0.5)$, while it has only one when   $(b,c,\rho)=(0.5,0.5,0.01)$. Nevertheless, in both cases, noise shifts the modes to cleaner states of the lake. In Figure \ref{change}, we illustrate the changes  induced to the value function by small changes in the rate of recycling, $r$. In particular, we numerically approximate the value functions $V$ that correspond to the rate of recycling  $r(x)=\frac{1}{2}(\tanh(a(x-3))+\tanh(3a))$ for various values of the parameter $a$, as well as the step function $\mathbf{1}\{x>3\}$ in  the stochastic $(\sigma = 0.1)$ case.

\section*{Acknowledgement} The authors are grateful to Emmanuil Georgoulis for his suggestions on the implementation of the numerical algorithm. 
\section*{Funding} This work has been supported by the Hellenic Foundation for Research and Innovation (H.F.R.I.) under the “First Call for H.F.R.I. Research Projects to support Faculty members and Researchers and the procurement of high-cost research equipment grant,” project HFRI-FM17-1034 (SCALINCS).

\vskip.05in
			
\bibliographystyle{apacite}
\bibliography{reference}
					
		\end{document}